\RequirePackage{fix-cm}
\documentclass[smallextended]{svjour3}       
\smartqed  
\usepackage{graphicx}
\usepackage{amsmath}
\begin{document}

\title{Simple algorithm for correction of a finite matrix group known approximately
  \thanks{This work is supported by INRTU grant for scientific groups No. 05-FPK-19
    }
}



\author{Andrey S. Mysovsky
}


\institute{
  Andrey S. Mysovsky \at
  $^1$Irkutsk National Research Technical University, 83 Lermontov Street, 664074 Irkutsk, Russia \\
  $^2$A.P. Vinogradov Institute of Geochemistry SB RAS, 1a Favorsky Street, 664033 Irkutsk, Russia \\
  \email{mysovsky@gmail.com}           
}

\date{Received: date / Accepted: date}

\maketitle

\begin{abstract}

  In this paper we have considered a finite unitary matrix group with exact elements being unknown and only approximate elements available. Such a group becomes inconsistent with its own  multiplication table. We found simple correction formula for such group. When applied iteratively this formula gives fast convergence of the group elements and allows to perform the group reconstruction.

  Next we considered small unitary rotation of entire group which makes the group consistent with a set of additional conditions.   Again, iterative procedure based on this correction shows fast convergence. All algorithms developed in this paper were implemented in a Python library which is available as open source software.

\keywords{finite groups \and group representations  \and molecular symmetry \and crystalline symmetry}
\subclass{20C35 \and 20H15  \and 15A30}
\end{abstract}

\section{Introduction}
\label{intro}

This paper deals with a situation when there is a finite group of unitary matrices (or operators)  $G=\{G_i\}$ which we need to know, but instead of them only approximate matrices $G_i'\approx G_i$ are available.
Besides the obvious disadvantage of not being exact, the set of $G'_i$ matrices has another serious drawback - it becomes inconsistent with G group multiplication table, is not closed any more with respect to multiplication operation and therefore does not comprise a group. In this paper the following problems are considered and solved:

\begin{itemize}
\item Problem 1. Perform the multiplication table based group reconstruction, i.e. find the matrices $G_i$ that are close to approximate matrices $G'_i$ but do obey the group multiplication table;

\item Problem 2a. In addition to problem 1, the matrices $G_i$ should be as close as possible to a set of given matrices $Q_i$ which are not required to be unitary and non-degenerate;

\item Problem 2b. In addition to problem 1, the matrices $G_i$ should turn a set of vectors $a_{ij}$ into vectors $b_{ij}$, i.e. approximate equality $G_i a_{ij} \approx b_{ij}$ should be fulfilled with best possible accuracy.
\end{itemize}

The solutions of all three problems have been implemented as an open source Python library freely available on GitHub \cite{grpcorr}. 

The example application considered in this paper concerns molecular symmetry when the geometry of the molecule is for various possible reasons distorted from symmetric configuration. These reasons can include thermal movement, presence of external field, spontaneous symmetry break, or be rather computational than physical: not fully converged geometry optimisation. However, the suggested algorithm can find much wider application.

The obvious generalisation of molecular and crystalline symmetry groups are the space groups in higher dimensions. In material science such groups are widely used to explain the structure of quasicrystals and other quasiperiodic solids. After the amazing discovery of a metallic alloy with crystallographically forbidden symmetry in 1984~\cite{shechtman} it was briefly recognised~\cite{janssen,janot} that its structure can be viewed as a slice of multidimensional periodic lattice projected into 3D space. Since that unusual properties of quasiperiodic solids continue to attract great attention~\cite{stein,delodi,liu,hua}, including recent discovery of superconductivity in quasicrystals~\cite{kamiya}. An up-to-date review can be found in Ref.~\cite{qrev}. An algorithm capable of producing a group out of set of matrices can potentially facilitate the recognition of quasicrystalline symmetries.

The paper is organised as follows. First, in section \ref{sec:multab} the problem 1 is considered. Second, in section \ref{sec:lsf} problems 2a and 2b are formulated as least-squares problems and common solution is suggested for them. A real example concerning molecular symmetry is provided in section \ref{sec:molsymm} together with  efficiency and convergence consideration of suggested method. After that the conclusion follows.

\section{\label{sec:multab}Multiplication table based group reconstruction}

To start with, let us introduce some notations. The group multiplication table will be denoted $(ij)$, i.e. multiplication of $i$-th and $j$-th group elements produces the element number $(ij)$. The group size is denoted $N$ and the deviation of approximate group matrices from exact ones is $\delta G_i=G'_i-G_i$. 

Let us introduce also $F_{ij}=G'_i \cdot G'_j - G'_{(ij)}$  - the matrices characterising the multiplication table violation. Up to the first order of magnitude over $\delta G_i$ they can be expressed as $F_{ij}=G_i\delta G_j+\delta G_iG_j-\delta G_{(ij)}$. Then we calculate

\begin{eqnarray}
  \label{fij_sum_1}
  \begin{aligned}
  \sum_j F_{ij}G_j^{-1} = N\delta G_i+G_i\sum_j \delta G_jG_j^{-1} \\
  -\sum_j \delta G_{(ij)}G_j^{-1}G_i^{-1}G_i = 
  N\delta G_i+[G_i,D_1]
  \end{aligned}
\end{eqnarray}

and

\begin{eqnarray}
  \label{fij_sum_2}
  \begin{aligned}
  \sum_j G_j^{-1}F_{ji} = N\delta G_i+\sum_j\delta G_jG_j^{-1}G_i \\
  -\sum_j G_iG_i^{-1}G_j^{-1}\delta G_{(ji)} = 
  N\delta G_i+[D_2,G_i]
  \end{aligned}
\end{eqnarray}

where $D_1=\sum_j\delta G_jG_j^{-1}$, $D_2=\sum_jG_j^{-1}\delta G_j$ and $[A,B]=AB-BA$ stands for the commutator. Taking sum of (\ref{fij_sum_1}) and (\ref{fij_sum_2}) allows us to get explicit expression for the deviations:

\begin{equation}
  \label{di_expression}
 2N\delta G_i = \sum_j(F_{ij}G_j^{-1}+G_j^{-1}F_{ji})+[G_i,\sum_j[G_j^{-1},\delta G_j]]
\end{equation}

By no means can the second term of this equation be expressed through the known values. In fact, it can be safely set to zero and there is a good reason for that.
 Indeed, consider the unitary rotation of entire group $G'_i = U^{-1}G_iU$ with a unitary matrix $U$ being the same for all $G_i$. Such a rotation does not violate the group multiplication table.  Up to the first order of magnitude small unitary rotation can be represented as $U=1+R$ with small anti-hermitian $R$. Thus, with the same accuracy

\begin{equation}
\label{grot}
\delta G_i = (1+R)^{+}G_i(1+R)  - G_i = [G_i,R]
\end{equation}

Here $'+'$ superscript stands for hermitian conjugate. The deviations of the form (\ref{grot}) are not distinguishable for multiplication table based group reconstruction and should be either set to zero or determined from some other sources.

Second term of (\ref{di_expression}) has exactly this structure. Indeed, since both $G_i$ and $G'_i$ are unitary, $(G'_i)^{+}G'_i-G_i^{+}G_i = \delta G_i^{+}G_i+G_i^{+}\delta G_i=0$ leading to the expression $\delta G_i^{+}=-G_i^{-1} \delta G_iG_i^{-1}$. It is trivial then to show that $[G_j^{-1},\delta G_j]$ commutators forming second term of (\ref{di_expression}) are anti-hermitian.

Substituting  $G_i$ (which are unknown) with $G'_i$ (which are known) in the expression (\ref{di_expression}) leads to the error only in the second order of magnitude over $\delta G_i$. We have arrived to the final first-order of magnitude expression for the deviations that solves Problem 1: 

\begin{equation}
  \label{multab_correction}
  \delta G_i = \frac{1}{2N}\sum\limits_j \left(F_{ij}(G'_j)^{-1}+(G'_j)^{-1}F_{ji}\right)
\end{equation}

The computational procedure based on this expression is simple and straightforward. It is applied iteratively until $G_i$ matrices are converged. Each iteration includes:

\begin{enumerate}
  
\item Re-unitarisation of $G_i$ matrices to enforce their unitarity. Re-unitarisation is done according to the formula $G_i \rightarrow (G_iG_i^{+})^{-1/2}G_i$;

\item Construction of multiplication table violation matrices $F_{ij}$;
  
\item Calculation of error value characterising the violation of multiplication table

  \begin{equation}
    \label{S_M}
    S_M=\sqrt{\sum_{ij} ||F_{ij}||^2}
  \end{equation}
  
  If the error is smaller than the convergence threshold value $\varepsilon$ then iterations are terminated and $G_i$ matrices are considered successfully converged;

\item Application of formula (\ref{multab_correction}) and subtraction of $\delta G_i$ deviations from $G_i$ matrices.
  
\end{enumerate}

The algorithm turns out to be quadratically convergent and converges very fast. The numerical consideration is provided below in section \ref{sec:molsymm}.

\section{\label{sec:lsf}Group correction with simultaneous least squares fit}

The  above described procedure allows one to find the $G_i$ matrices obeying the group multiplication table. However, these matrices are still not uniquely defined. They are defined only up to arbitrary unitary rotation of entire group, and this remaining degree of freedom can be used to solve Problems 2a and 2b. First, we formulate them as least-squares problems:

\begin{itemize}

  \item $\sum_i ||G_i - Q_i||^2 = min$ (problem 2a);
  \item $\sum_{ij}||G_i a_{ij} - b_{ij}||^2 = min$ (problem 2b);
\end{itemize}

Here $||A||^2=\sum_{\alpha \beta} |A_{\alpha \beta}|^2 = Tr (A^{+}A)$ is the matrix norm squared. Now we consider unitary rotation of the group $G_i \rightarrow e^{-R}G_i e^{R}$ expressed through small anti-hermitian matrix $R$. Anti-hermiticity condition $R+R^{+}=0$ should be used as a constraint when minimising the sum of squares for problems 2a and 2b. Our goal is to find $R$ matrix minimising the corresponding error function (sum of squares).



Now we are ready to write down the value which is to be minimised in order to solve Problem 2a:

\begin{eqnarray}
  \label{serr}
    S = \sum_i Tr \left[(e^{-R}G_i e^{R} - Q_i)^{+}(e^{-R}G_i e^{R} - Q_i)\right] 
    - Tr \left[\lambda (R+R^{+})\right]
\end{eqnarray}

Here $\lambda$ is the Lagrange multipliers matrix to account for the constraint. In order to obtain linear equation for $R$, this expression has to be expanded up to second order over $R$ and subsequently differentiated:

\begin{eqnarray}
  \label{dsdr}
  \begin{aligned}
    \frac{\partial S}{\partial R} = \sum\limits_i \left \{  2(R+R^{+}) -
    G_i^{+}(R+R^{+})G_i  
    - G_i(R+R^{+})G_i^{+}  + [G_i,Q_i^{+}] \phantom{\frac{1}{2}} \right. \\
    \left. - \frac{1}{2}\left( R(G_iQ_i^{+}+Q_i^{+}G_i) + (G_iQ_i^{+}+Q_i^{+}G_i)R
    - 2 G_i R Q_i^{+} - 2 Q_i^{+} R G_i\right) \right\} -\lambda =0
  \end{aligned}
\end{eqnarray}

Similarly, one can proceed with the differentiation over $R^{+}$:

\begin{eqnarray}
  \label{dsdrcross}
  \begin{aligned}
    \frac{\partial S}{\partial R^{+}} = \sum\limits_i \left \{  2(R+R^{+}) -
    G_i^{+}(R+R^{+})G_i  
    - G_i(R+R^{+})G_i^{+}  - [G_i^{+},Q_i] \phantom{\frac{1}{2}} \right. \\
    \left. - \frac{1}{2}\left( R^{+}(G_i^{+}Q_i+Q_iG_i^{+}) + (G_i^{+}Q_i+Q_iG_i^{+}) R^{+}
    - 2 G_i^{+} R^{+} Q_i - 2 Q_i R^{+} G_i^{+}\right) \right\} -\lambda =0 
  \end{aligned}
\end{eqnarray}

Subtraction of equations (\ref{dsdr}) and (\ref{dsdrcross}) allows one to obtain simpler equation, and, what is more important, to get rid of Lagrange multipliers $\lambda$. After that it is safe to substitute $R^{+}\rightarrow -R$ everywhere and get the final equation for R matrix:

\begin{eqnarray}
  \label{req}
  \sum_i \left\{ RH_i + H_iR - (G_iRQ_i^{+}+Q_i^{+}RG_i + G_i^{+}RQ_i + Q_iRG_i^{+})\right\}
  \nonumber \\
  = \sum_i [G_i,Q_i^{+}] + [G_i^{+},Q_i]
\end{eqnarray}

with $H_i$ denoting the following structure:

\begin{equation}
  \label{defh}
  H_i = \frac{1}{2}\left( G_iQ_i^{+} + Q_i^{+}G_i +Q_iG_i^{+}+G_i^{+}Q_i \right)
\end{equation}

The matrix equation (\ref{req}) 
falls into the class of Generalised Sylvester equations \cite{hautus}. There is a straightforward but possibly not very efficient solution of such equation. It can be seen as a system of linear equations for elements of matrix R. To re-formulate it in this manner, we need the matrix vectorisation operator $\text{vec}\,R$ that puts all columns of R in the same vector one after another and the Kronecker product operation

$$
A\otimes B = \begin{pmatrix}
  A_{11}B & \dots & A_{1n}B\\
  \hdotsfor{3} \\
  A_{n1}B & \dots & A_{nn}B\\
  \end{pmatrix}
$$

Now (\ref{req}) is equivalent to the linear equations system $L\cdot \text{vec}\, R = \text{vec}\,P$ with the right-hand side denoted

\begin{equation}
  \label{rightp}
  P=\sum_i [G_i,Q_i^{+}] + [G_i^{+},Q_i] 
\end{equation}

and $L$ being the ``supermatrix'':

\begin{equation}
  \label{supmtr}
  L = \sum_i \left\{ H_i^{T} \otimes I_n + I_n \otimes H_i
  -(Q_i^{*}\otimes G_i + G_i^{T}\otimes Q_i^{+} + Q_i^{T}\otimes G_i^{+} + G_i^{*}\otimes Q_i) \right\}
\end{equation}

Here $'T'$ and $'*'$ superscripts stand for the transpose and complex conjugate of matrices, respectively, while $I_n$ is the identity matrix of dimension $n$. Formulae (\ref{defh},\ref{rightp},\ref{supmtr}) together with the expression $\text{vec}\,R = L^{-1}\cdot \text{vec}\,P$ provide the formal solution of Problem 2a.

However, this solution can be complicated by the presence of zero eigenvalues of $L$ supermatrix. Such zero eigenvalues correspond to eigenmatrices which, being substituted instead of $R$ into left-hand side of (\ref{req}), turn it to zero. There is always at least one such eigenmatrix with zero eigenvalue - it is the identity matrix. However, it turns out that all such eigenmatrices are hermitian and can be easily ruled out when constructing the solution of (\ref{req}).

Another serious drawback of the suggested solution is its scalability. Indeed, if the matrices $G_i$ have the dimension $n\times n$, then $L$ supermatrix has already $n^2\times n^2$. Taking into account that commonly used eigensolvers for dense matrices scale cubically, that gives us total scalability of $n^6$. For large $n$ it can be computationally very expensive and even prohibitive.

However, a great simplification of this solution can be done if it can supposed that $G_i \approx Q_i$. Denoting $\delta Q_i = Q_i - G_i$ one can see that right-hand side of (\ref{req}) depends on these quantities linearly: $P = \sum_i [G_i,\delta Q_i^{+}] + [G_i^{+}, \delta Q_i] $. After that substituting $Q_i$ with $G_i$ on the left-hand side of this equation creates an error only of second order of magnitude over $\delta Q_i$.

Then immediately each $H_i=2$ because of the unitarity of $G_i$ matrices.
Using once again the fact that $G_i^{+}=G_i^{-1}$ and reordering the summation in some terms at the left-hand side of equation (\ref{req}) gives us:

$$
4NR - 4 \sum_i G_iRG_i^{+} = P
$$

Consider the eigenmatrices of left-hand side of this equation. Such an eigenmatrix should obey

$$
4NR - 4 \sum_i G_iRG_i^{+} = \alpha R = (4N - 4 \mu) R
$$

with eigenvalue $\alpha$ which we represent as $\alpha=4(N-\mu)$. Then $\mu$ is an eigenvalue for another matrix transformation: $\sum_i G_i R G_i^{+} = \mu R$. Multiplying by $G_j$ on the left side and by $G_j^{+}$ on the right side and reordering the summation leads us to the conclusion that $\mu R = \mu G_j R G_j^{+}$ for any $j$, in other words, that the matrix $\mu R$ must commute with the entire group $G$.

The matrices of the group $G$ themselves comprise some representation of this group. It can be either irreducible or reducible, in that latter case they can be decomposed into the direct sum or irreducible representations: $G_i=G_i^{(1)}\oplus G_i^{(2)}\oplus\dots\oplus G_i^{(k)} $. In other words, there is a unitary basis set transformation U which makes all $G_i$ matrices block-diagonal with the blocks being the matrices of corresponding irreducible representation. According to the Schur's lemma any matrix commuting with irreducible representation matrices must be identity matrix multiplied by a scalar.

Thus, we arrived to the conclusion that $\mu R = U^{+}( \mu_1 I_{n_1} \oplus  \mu_2 I_{n_2} \oplus \dots \oplus \mu_k I_{n_k} )U$. If this construction is non-zero, then it is necessarily hermitian. For any anti-hermitian matrix $R$ it means that $\mu=0$, and we can get simple explicit expression for $R$:

\begin{equation}
  \label{simple_req}
  R = \frac{1}{4N}\sum_i \left( [G_i,Q_i^{+}] + [G_i^{+},Q_i] \right)
\end{equation}

This formula, although being approximate as the others in the present paper, is correct up to first order over $\delta Q_i$ and is free from the drawbacks of previous supermatrix-based solution.

It remains only to show that Problem 2b can be reduced to Problem 2a.  Error function (sum of squares) for Problem 2b is

\begin{eqnarray*}
    S = \sum_{ij} <G_ia_{ij}-b_{ij}|G_ia_{ij}-b_{ij} > - Tr\left[\lambda(R+R^{+})\right] \\ =
    \sum_i Tr \left\{ G_i \left(\sum_j |a_{ij}><a_{ij}|\right)G_i^{+}
    - G_i \left(\sum_j |a_{ij}><b_{ij}| \right) \right. \\
    \left. - \left(\sum_j |b_{ij}><a_{ij}|\right)G_i^{+}  
    + \left(\sum_j |b_{ij}><b_{ij}| \right)\right\} - Tr\left[\lambda(R+R^{+})\right]
\end{eqnarray*}

If we set

\begin{equation}
  \label{Q_ab}
  Q_i = \sum_j |b_{ij}><a_{ij}|
\end{equation}
  
then second and third terms under the sum over $i$ become identical to those in (\ref{serr}), while the first and the fourth do not, but they make no contribution into equation
(\ref{req}). This means that Problem 2b can be solved identically to Problem 2a with such a choice of $Q_i$ matrices. 

The algorithm solving these two problems is, again, iterative. Each iteration consists of the following steps:

\begin{enumerate}
  
\item Re-unitarisation of $G_i$ matrices;

\item Calculation of multiplication table violation error $S_M$ given by (\ref{S_M}). If this value exceeds the threshold $\varepsilon_M$, one step of multiplication table based reconstruction is performed and the convergence flag is set to ``False'';
  
\item Calculation of $P$ matrix via (\ref{rightp});

\item Calculation of $R$ by either supermatrix-based approach ($\text{algo}=1$) or by simplified formula (\ref{simple_req}) ($\text{algo}=2$). The algo parameter is the choice of the user.

\item Rotation of the group as $G_i \rightarrow e^{-R} G_i e^R$;

\item Calculation of error value characterising the difference between $G_i$ and $Q_i$ matrices

  \begin{equation}
    \label{S_Q}
    S_Q=\sum_i ||G_i-Q_i||^2\text{ (for problem 2a)}
  \end{equation}

or
  
  \begin{equation}
    \label{S_ab}
    S_Q=\sum_{ij} ||G_ia_{ij}-b_{ij}||^2\text{ (for  problem 2b)}
  \end{equation}

  and its change on current iteration $\Delta S_Q = S_Q(prev) - S_Q(current)$;

\item Comparison of $\Delta S_Q$ with the threshold value $\varepsilon_Q$ and comparison of $||R||$ with $\varepsilon_R$. If one of these values exceeds convergence threshold, the convergence flag is set to ``False'';

\item Checking the convergence flag. If it is ``True'', the fit is successfully converged.
  
\end{enumerate}

\section{\label{sec:molsymm}Application to molecular symmetry}

\begin{figure*}
  \begin{tabular}{l l}
    a) & b) \\
    \includegraphics[scale=0.5]{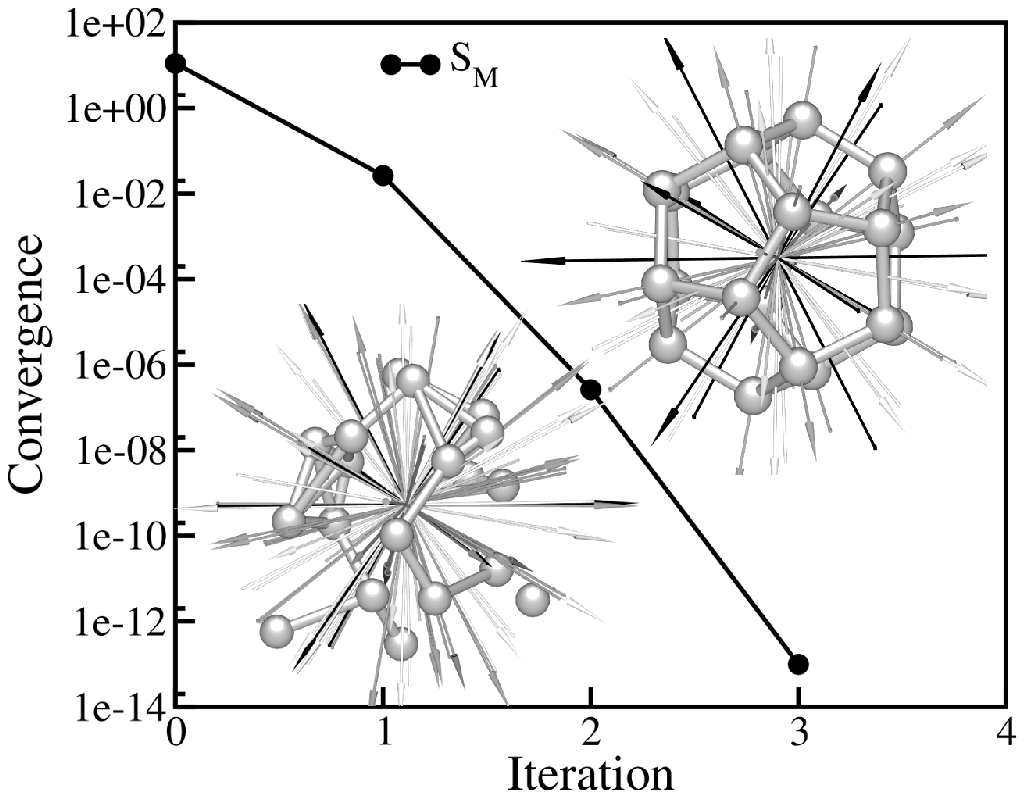} &
    \includegraphics[scale=0.5]{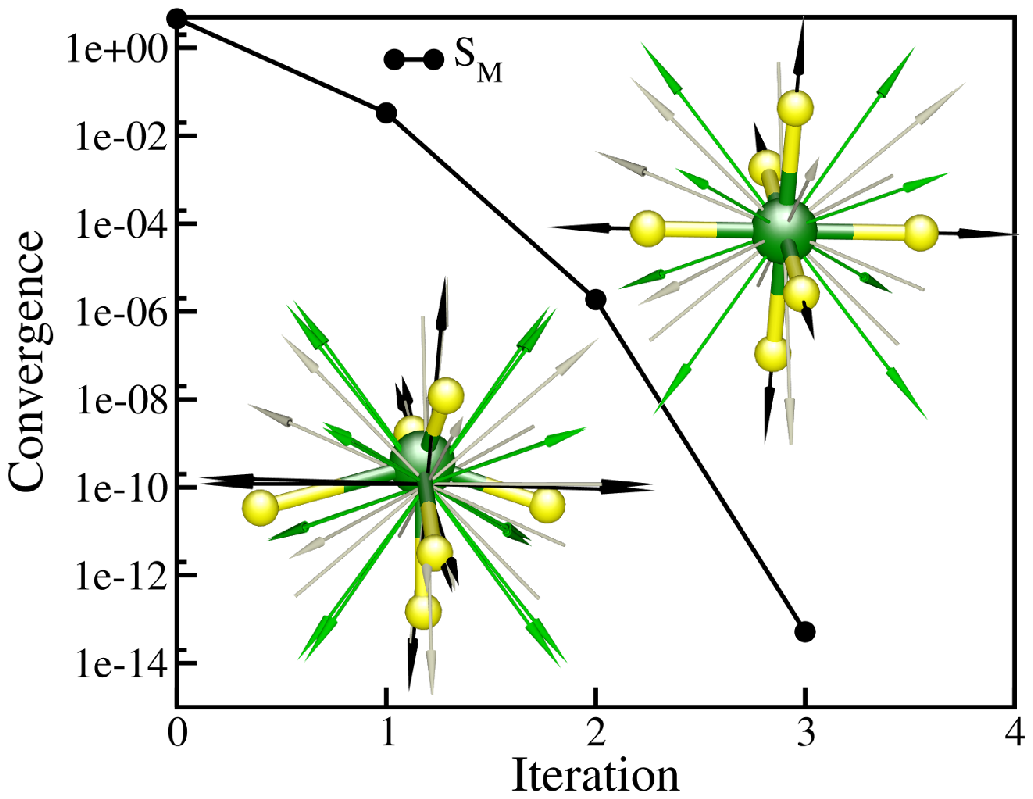} \\
  \end{tabular}
  \caption{\label{fig:epsart} Reconstruction of randomly distorted a) C$_{20}$ fullerene and its symmetry group I$_h$; b)  SF$_6$ molecule and its symmetry group O$_h$.}
\end{figure*}


Application of the reported methods to the symmetry of molecules, nanostuctures and crystals has been in fact the main driving force behind this work, although their applicability turned out to be much wider. In this section we consider how the suggested algorithms work for the groups of molecular symmetry. Through this section we use the Sh\"{o}nflies notation for point groups~\cite{flurry}. 


First, the multiplication table based group reconstruction (the solution of problem 1) has been tested. Two molecules, namely, C$_{20}$ fullerene and sulphur hexafluoride were selected for this test. Their symmetry groups are I$_h$ and O$_h$, respectively. The molecules were constructed as simple list of Cartesian coordinates of each atom, the length unit for the coordinates was \AA ngstr\"{o}m. 

Multiplication table for the symmetry group of each molecule has been constructed as follows. All possible permutations of atoms that do not change interatomic distances and bond angles comprise the permutation group isomorphic to the point symmetry group of the molecule. This permutations group was used to build the multiplication table. After that geometries of both molecules were distorted by adding random displacements to each coordinate. Each displacement value for each coordinate ranged uniformly from -0.5 to 0.5~\AA. Next, the group matrices of initial approximation were constructed via formula (\ref{Q_ab}) and re-unitarised. Finally, the multiplication table based group reconstruction has been applied. The driver routine implemented in the repository~\cite{grpcorr} is called {\it multab\_group\_correction}.

The results are shown on fig. 1. It requires only 4 iterations to achieve the convergence of $10^{-12}$ for the multiplication table violation error $S_M$. Note also that the convergence becomes faster with each next iteration due to the quadratucally convergent character of the suggested algorithm.

The insets on both figs. 1a and 1b show the randomly distorted molecule before reconstruction and symmetrised molecule after the group reconstruction. The rotation axis of the symmetry group are depicted with arrows. The highest order axis (5-th for C$_{20}$ and 4-th for SF$_6$) are shown black while the others have different colours. Observe
 the bundles of axis  before the reconstruction turns them into just one axes. Due to inconsistency of approximate group some rotations which should have exactly the same axes have close but different ones. Reconstruction algorithm successfully corrects both rotation axis and angles.

\begin{figure}
\includegraphics[scale=0.5]{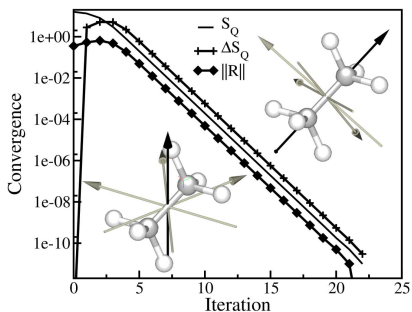}
\caption{\label{fig:epsart} Least-squares rotation fit  for D$_{3d}$ symmetry group of C$_2$H$_6$ molecule.}
\end{figure}

Next, the least-squares fit algorithm has been tested. As an example fig. 2 shows the ethane molecule. In the beginning the molecule had its full symmetry  D$_{3d}$ which was correctly determined similarly to previous example. Then C$_2$H$_6$ molecule was rotated around a random axis by random angle (on fig. 2 it was 86.47$^{\circ}$). After that the algorithm solving Problem 2b has been applied (the corresponding driver routine in the repository~\cite{grpcorr} is called {\it abfit\_group\_correction}). The convergence on fig. 2 is characterised by previously defined quantities $S_Q$ calculated according to (\ref{S_ab}), its change and the rotation norm $||R||$. As one can see, all three quantities decay exponentially during the fit and after it is finished the symmetry group axes are oriented properly with respect to the molecule.

\begin{figure}
\includegraphics[scale=0.5]{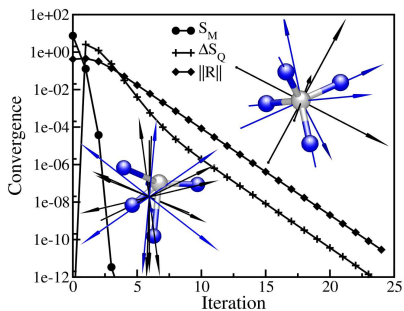}
\caption{\label{fig:epsart} Reconstruction with simultaneous least-squares rotation fit  of T$_d$ group for CH$_4$ molecule.}
\end{figure}

Finally, we have tested the ability of suggested algorithms to manage both types of error simultaneously, i.e. to correct the group multiplication table and at the same time find the unitary rotation. This is illustrated by fig. 3. CH$_4$ molecule has been randomly deformed (each coordinate received random displacement from -0.5 to 0.5~\AA). Then initial  group matrices were constructed by formula (\ref{Q_ab}) and re-unitarised. After that the molecule was randomly rotated (by 56.93$^{\circ}$ on fig. 3). When all this preparations were finished, the algorithm for Problem 2b has been applied.

The convergence was fastest for $S_M$ value similar to previous examples. After only 4 iterations it became virtually zero. The values of $\Delta S_Q$ and $||R||$ showed exponential decay and converged to threshold value during 24 iterations. In this example the value of $S_Q$ does not converge to zero because the molecule is not only rotated but also deformed. Instead it converges to its minimal value (8.316 for the example shown on fig. 3).

It is interesting to note that algorithm for Problem 2b always finds the correct rotation only when the rotation angle is less then 90$^{\circ}$. For larger angles the convergence is never achieved. As the test calculations described in this section were repeated several hundreds of times this can be stated with confidence. Another interesting observation is that there is almost no difference in the convergence between the ``supermatrix'' approach and simplified formula~(\ref{simple_req}) for the problems 2a and 2b.

\section{\label{sec:concl}Conclusion}

In this paper we have considered a finite unitary matrix group with exact elements being unknown and only approximate elements available. The first problem of such a group is the inconsistency with its own multiplication table which we called ``Problem 1''. Another distinct inconsistency type arises when it is known how the group elements should act on certain set of vectors, but they do not act as they should. We formulated this as ``Problem 2a'' and ``Problem 2b''.

The solutions of all three considered problems were found analytically as a first order correction to the group elements. Those solutions can be used in iterative procedure until the group matrices are converged, and such iterative procedures were implemented in the GitHub repository~\cite{grpcorr}.
The present Python implementation does not pretend to be highly computationally efficient. It was programmed rather to test the suggested methods and explore their applicability. However, as it is published as open-source software, any researcher can feel free to use the code or to rewrite it in more efficient manner.

The numerical test performed for the point symmetry groups had shown very good convergence of suggested methods. The multiplication table based reconstruction has quadratic convergence, while for the least-squares fit of unitary rotation it is linear. Although initially we considered the situation when the deviations of approximate group matrices from exact ones are small, the suggested algorithms seem to successfully handle the situations when they are not.


\begin{thebibliography}{}
%
%
\bibitem{grpcorr}
  A.S. Mysovsky, Python library for reconstruction of a finite matrix group known approximately, GitHub repository, \texttt{https://github.com/mysovsky/grpcorr} (2019)
  
\bibitem{shechtman}
D. Shechtman, I. Blech, D. Gratias, and J. W. Cahn, Metallic Phase with Long-Range Orientational Order and No Translational Symmetry, 
Phys. Rev. Lett. 53, 1951 (1984)

\bibitem{janssen}
  T. Janssen, Crystallography of Quasi-Crystals, Acta Cryst. A42, 261-271 (1986)
  
\bibitem{janot}
  Ch. Janot and J.M. Dubois, Quasicrystals, J. Phys. F: Met. Phys., 2303-2343, 18 (1988)
  
\bibitem{stein}
  P. J. Steinhardt and H.-C. Jeong,A simpler approach to Penrose tiling with implications for quasicrystal formation,  
  Nature 382, 431–433 (1996) 
  
\bibitem{delodi}
S. Deloudi and W. Steurer, Higher-dimensional crystallography of N-fold
quasiperiodic tilings, Acta Cryst. A68, 266-277 (2012)

\bibitem{liu}
  J.F. Liu, Z.Q. Yang, H.Q. Ye,
  Solid-state formation of icosahedral quasicrystals at Zn3Mg3Y2/Mg interfaces in a Mg–Zn–Y alloy,
  J. Alloy Compd. 650, 65-69 (2015)

\bibitem{hua}
Hua Li, Haikun Ma, Longgang Hou, Mingxing Guo, Zhanbing He,
Shield-like tile and its application to the decagonal quasicrystal-related structures in Al-Cr-Fe-Si alloys,
J. Alloy Compd. 701, 494-498 (2017)

\bibitem{kamiya}
K. Kamiya, T. Takeuchi, N. Kabeya, N. Wada, T. Ishimasa, A. Ochiai, K. Deguchi, K. Imura and N. K. Sato,
Discovery of superconductivity in quasicrystal, 
Nat. Commun. 9, 154 (2018)

\bibitem{qrev}
W. Steurer, Quasicrystals: What do we know? What do we
want to know? What can we know?, Acta Cryst. A74, 1–11 (2018)

\bibitem{hautus}
M.L.J. Hautus, Operator substitution, Linear Algebra Appl. 205-206,  713-739 (1994)
  
\bibitem{flurry}
R.L. Flurry, Symmetry Groups : Theory and Chemical Applications. Prentice-Hall Inc, New Jersey (1980)
  
\end{thebibliography}
\end{document}